\documentclass[11pt]{article}
\usepackage{latexsym,amsmath,amsthm,amssymb,verbatim,ifthen,amssymb,graphicx,color,mathrsfs}

\usepackage[english]{babel}
\usepackage[utf8]{inputenc}

\usepackage{tikz-cd}
\usetikzlibrary{arrows.meta,decorations.markings}
\tikzset{
	disk vertex/.style={circle, fill, inner sep=1.5pt},
	disk edge/.style={thick, postaction={decorate},
		decoration={markings, mark=at position 0.5 with
			{\arrow{Stealth[length=2.4mm, width=2mm]}}}},
	symmetry/.style={very thick},
	pics/celldisk/.style={code={
		\fill[gray!20] (0,0) circle (1.5);
		\draw[disk edge] (180:1.5) arc[start angle=180, end angle=0, radius=1.5];
		\draw[disk edge] (0:1.5) arc[start angle=0, end angle=-180, radius=1.5];
		\node[disk vertex, label=left:$b$] at (180:1.5) {};
		\node[disk vertex, label=right:$a$] at (0:1.5) {};
		\node[above=2pt] at (90:1.5) {$x$};
		\node[below=2pt] at (-90:1.5) {$y$};
	}},
}
\usepackage{mathtools}  

\usepackage[all]{xy}
\usepackage{color}

\usepackage{hyperref}
\hypersetup{
	colorlinks=true,
	allcolors=blue,
}  		
\usepackage{anysize}
\usepackage{cancel}									
\marginsize{2.5cm}{2.5cm}{2cm}{2cm}
\usepackage{enumerate}

\SelectTips{cm}{}

\newtheorem{theorem}{Theorem}[section]

 \newtheorem{lemma}[theorem]{Lemma}
 \newtheorem{proposition}[theorem]{Proposition}
 \theoremstyle{definition}
 
 \theoremstyle{remark}
 \newtheorem{remark}[theorem]{Remark}
 \newtheorem{example}[theorem]{Example}
 \numberwithin{equation}{section}

\DeclareRobustCommand{\stirling}{\genfrac\{\}{0pt}{}}

\newcommand{\CompLie}{{\widehat{\mathbb{L}}}}

\DeclareMathOperator{\ad}{ad}

\begin{document}

\title{Identities for generalized Bernoulli numbers\\ from a Lie model of the 2-disk}

\author{Urtzi Buijs, José Manuel Moreno Fernández}
\date{}

\maketitle

\begin{abstract}
\emergencystretch=1em
We define generalized Bernoulli numbers as appropriate coefficients of
the Baker--Campbell--Hausdorff (BCH) product formula, extending the
classical sequence of rational numbers to a 2-dimensional setting.
We give an identity involving these numbers and the classical
Bernoulli numbers which is equivalent to the existence of a particular
complete differential graded Lie algebra modelling a cell
decomposition of the 2-disk and based on the Lawrence--Sullivan model
of the interval. From this identity we derive explicit quadratic
relations.

The existence of Lie models for the symmetries of the interval and the
disk is related to certain properties of the BCH product and to some
notable features of the classical Bernoulli numbers: the two choices
for $B_1$, the vanishing of $B_n$ for odd $n\geq 3$, and the equality
$B_4=B_8$.
\end{abstract}

\section{Introduction} 
\label{sec: Introduction}

The Bernoulli numbers are arranged in a sequence $B_n$ of rational numbers appearing all across mathematics, 
both in classical settings and in contemporary research.
Several references deal with many of the well-known properties and identities of this sequence, 
see for instance \cite{Con96,BerAndZeta}.
These properties and relationships have been studied
using tools in several fields of mathematics,
such as number theory, combinatorics, or analysis.

\smallskip

In this work, 
we emphasize that certain relationships among the Bernoulli numbers and a certain generalization of these
can be derived by using  
\emph{topological} and 
\emph{discrete differential geometry} arguments. 
This program was already started in \cite{BCM}.
There, it was proven that the equation $$\partial^2 = 0$$
arising from the differential on a very special complete differential graded Lie algebra $\left(\mathcal L_1,\partial\right)$ due to Lawrence--Sullivan \cite{LawSul} is equivalent to a certain family of quadratic identities among the Bernoulli numbers 
of the sort 
$$\sum_{\substack{k = 0 \\ k \textrm{ even}}}^n \lambda_k B_k B_{n-k} = 0.$$
Euler-type relations among the Bernoulli numbers had also been
obtained by Parent and Tanr\'e as a by-product of their proof that the
two Lie models of the interval constructed by Lawrence and Sullivan
coincide \cite[Section 4]{PT}.

\medskip

In this paper,
we take a step further and study the identities among (generalized)
Bernoulli numbers encoded in the complete dg Lie algebra built on the
cell decomposition of the $2$-disk with two vertices, two edges and
one $2$-cell \cite{GGL}.
The existence of the differential of this complete dg Lie algebra,
that is, the equation $\partial^2\varepsilon=0$ on the $2$-cell, is
equivalent to an identity involving the generalized Bernoulli numbers
and the classical ones (Remark \ref{rem:equivalence}). From this
identity we derive explicit quadratic relations, which are collected
in Theorem \ref{teo: mainintro}, the main result of this paper.
Since the explicit formulas are convoluted and require setting up a
bit of notation, we refer the reader to
Section \ref{sec: generalized Bernoulli and 2-disk}
instead of reproducing the result here.
In turn, the equation $\partial^2\varepsilon=0$ is equivalent to
$\partial(x*y)=[x*y,a]$, which expresses the gluing of two
Lawrence--Sullivan intervals along their end points
(Theorem \ref{subdivision}), so the $2$-cell plays a formal role
here.

\smallskip

The Bernoulli numbers can be defined, 
up to a division by a factorial, 
as the coefficients appearing in the expansion of the 
exponential generating function,
$$\frac{x}{e^x-1} = \sum_{n\geq 0} \frac{B_n}{n!}x^n.$$
In connection with the Baker--Campbell--Hausdorff
product,
$$x* y=x+y+\frac{1}{2}[x,y]+\frac{1}{12}[x,[x,y]]-\frac{1}{12}[y,[x,y]]+\cdots ,$$
it is well known that the Bernoulli numbers are recovered 
(modulo factorials)
as the coefficients of the Lie monomials containing the variable $x$ exactly once.	
Following this observation, 
let $$w = x^{h_1}y^{k_1} x^{h_2}y^{k_2}\cdots x^{h_n}y^{k_n}$$
be a word in $x$ and $y$, viewed in the complete tensor algebra on two variables, $\widehat T(x,y)$.
This associative algebra,
which we assume throughout to have rational coefficients,
is the universal enveloping algebra of the complete free Lie algebra on two variables, $\CompLie\left(x,y\right)$.
This is the natural algebraic structure that supports the Baker--Campbell--Hausdorff product formula.
We define the \emph{generalized Bernoulli number} 
associated to the word $w$ as the rational number 
$$\mathcal{B}_w = h_1!k_1!\cdots h_n!k_n!c_{\left( \left(h_1, k_1\right),\dots , \left(h_n, k_n\right)\right)},$$
where $c_{\left( \left(h_1, k_1\right),\dots , \left(h_n, k_n\right)\right)}\in \mathbb Q$ is the coefficient of the word 
$w$ in the expansion of $x*y$.
Our main result, 
Theorem \ref{teo: mainintro},
gives some explicit quadratic relationships among these numbers. 
We compute some explicit examples in Section \ref{sec: examples}.

\smallskip

Finally, in Section \ref{sec: symmetries} we relate geometric
properties of the interval and of the disk to some notable facts
about the Bernoulli numbers. Namely, we illustrate the following
facts.
{\settowidth{\leftmargini}{Fact 3.}\addtolength{\leftmargini}{\labelsep}
\begin{enumerate}
	\item[Fact 1.] The vanishing of the odd terms $B_{2n+1}$ of the Bernoulli sequence 
	except for the first term $B_1$, 
	corresponds to the symmetry of the unit interval.	
	\item[Fact 2.] The existence of two choices for the initial term $B_1$ of the sequence 
	(namely, $\frac{1}{2}$ or $-\frac{1}{2}$) 
	corresponds to the existence of two opposite orientations on the unit interval.
	\item[Fact 3.] The rotation of the disk by $\pi$, which
	interchanges the vertices and also the edges, implies, together
	with Fact 1, the equality $B_4=B_8$. This is the only coincidence
	between two non-zero Bernoulli numbers with different indices.
\end{enumerate}}

\medskip

\noindent {\bf Notation and conventions.}
In this work, 
$I$ denotes both the closed unit interval $[0,1]$ as a topological subspace of the real line, 
as well as the standard 1-dimensional simplicial set $\Delta^1_\bullet$.
It will be clear from the context which one is in use.
All algebraic structures are taken with coefficients over
the rational numbers.

\bigskip

\noindent 
\textbf{MSC 2020:} 55P62, 11B68, 17B01

\smallskip
\noindent 
\textbf{Key words and phrases:} 
Rational homotopy theory, Bernoulli numbers, complete Lie algebra, Baker--Campbell--Hausdorff formula

\smallskip

\noindent {\bf Acknowledgments:}
The second author has been partially supported by the MICINN grant PID2023-149804NB-I00,
the Junta de Andalucía grant ProyExcel-00827, 
and the Junta de Andalucía research group grupoPaidi-G-FEDER-FQM264.

\smallskip

\noindent {\bf Disclaimer:}
AI was used solely in the final stage of the process to proofread the final version.

\subsection{Preliminaries}
\label{Preliminaries}

In this section, 
we collect the necessary background on the Baker--Campbell--Hausdorff product (BCH product, henceforth) and
the properties of it that will be used,
and on complete dg Lie models.
The background on these topics can be found for instance in 
the monograph \cite{LieModels}.
The generalized Bernoulli numbers are defined in
Section \ref{sec: generalized Bernoulli and 2-disk}. Up to the
factorials of the exponents, they are the coefficients of the BCH
series in $\widehat T(x,y)$, which were determined by Goldberg
\cite[Theorem 1]{Goldberg}.

\subsubsection{The Bernoulli numbers}

The Bernoulli numbers $B_n$ can be defined in multiple ways.
For instance,
in terms of the generating function
$$\frac{x}{e^x-1} = \sum_{n\geq 0} \frac{B_n}{n!}x^n.$$
In Lie theory they appear, up to signs and factorials, as the
coefficients of $\ad_x/(1-e^{-\ad_x})$, the inverse of the operator
$(1-e^{-\ad_x})/\ad_x$ that gives the differential of the
exponential map at $x$.

\subsubsection{The Baker--Campbell--Hausdorff product}
\label{sec: BCH prerequisites}
The BCH (\emph{Baker--Campbell--Hausdorff}) product is the group multiplication $*$ that arises when integrating, in the most classical sense, a nilpotent Lie algebra $L$ into a Lie group $\exp(L)$. 
In fact, 
one can integrate \emph{complete}
Lie algebras over a characteristic zero field into Malcev complete groups, producing a well-known equivalence of categories,
see \cite[Chapter 8]{Fresse1} for a modern account.
In this section, we briefly summarize the details of this topic that will be relevant to achieve our results.

\medskip

Recall that any Lie algebra $L$ can be canonically filtered by its lower central series,
\[
L = F^1L \supseteq[L,L]= F^2L \supseteq [L,[L,L]] = F^3L  \supseteq \cdots.
\]
By definition, the Lie algebra $L$ is complete if the natural map 
\[
L \longrightarrow \widehat L
\]
from $L$ to its completion with respect to the lower central series, $ \widehat L = \varprojlim L/F^nL$, is an isomorphism.
When $L$ is filtered by the lower central series, 
its universal enveloping algebra $UL$ acquires a corresponding natural filtration.
For a complete Lie algebra $L$, we write $UL$ for the completion of
its universal enveloping algebra with respect to this filtration. For
the complete free Lie algebra on $x$ and $y$, this is the complete
tensor algebra $\widehat T(x,y)$.
Given a complete Lie algebra $L$,
one defines the exponential and logarithm maps 
\begin{center}
	\begin{tikzcd}
		\exp : \widehat{\bar I}  \arrow[r, shift left] & 1 + \widehat{\bar I} : \log, \arrow[l, shift left]
	\end{tikzcd}
\end{center}
where $\widehat{\bar I}$ is the completion of the augmentation ideal $\bar I$ of the universal enveloping algebra $UL$ of $L$ with respect to its iterated powers.
The maps above
are inverse bijections, 
and are 
given by the classical power series expansions
$$\exp(x) = \sum_{n\geq 0} \frac{x^n}{n!} \qquad \textrm{and} \qquad \log\left(1+x\right) = \sum_{n \geq 1} \left(-1\right)^{n+1} \frac{x^n}{n}.$$
The BCH product $*$ on two elements  $x,y \in L$ is then defined as
$$x* y=\log\left(\exp(x)\exp(y)\right).$$
A well-known explicit formula due to Dynkin  in terms of nested commutators \cite{Dynkin} begins by 
$$x* y=x+y+\frac{1}{2}[x,y]+\frac{1}{12}[x,[x,y]]-\frac{1}{12}[y,[x,y]]+\cdots .$$
This formula   lives naturally in the complete free Lie algebra $\CompLie\left(x,y\right)$ generated by $x,y$, 
whose universal enveloping algebra is the complete tensor algebra $\widehat T\left(x,y\right)$.
In this associative algebra, 
whose product we denote by  juxtaposition,
expanding $\log(\exp(x)\exp(y))$ gives
$$x* y = \sum_{n\geq 1} z_n(x,y),$$
where 
$$z_n(x,y) = \sum_{k=1}^n\frac{(-1)^{k-1}}{k}\sum \frac{x^{p_1}y^{q_1}\cdots x^{p_k}y^{q_k}}{p_1!q_1!\cdots p_k!q_k!}.$$
The second sum above is taken over all possible partitions of $n$ by non-negative integers $(p_1,...,p_k,\allowbreak q_1,...,q_k)$ satisfying $p_i + q_i \geq 1$ for all $i$.
Thus, $p_1+q_1+\cdots +p_k+q_k=n.$

\smallskip

For any $x\in L$,
denote by $\ad_x = [x,-]$
the adjoint action, or inner derivation, determined by $x$.
We use the same notation for this map both in $L$ and $UL$.

Later on, 
we will need to consider complete \emph{differential graded} (dg) Lie algebras.
A dg Lie algebra is a $\mathbb{Z}$-graded chain complex
$L = \left(L,\partial\right)$
together with a degree zero graded antisymmetric Lie bracket $[-,-]$ satisfying the graded Jacobi identity, 
and which is compatible with the differential via the graded Leibniz rule, $$\partial[x,y] = [\partial x,y] + (-1)^{|x|}[x,\partial y]$$ 
for all homogeneous $x,y\in L$. 
Here, $|x|$ is the degree of $x$.
These algebras will always be assumed to be filtered by the lower central series,
which is actually a filtration of chain complexes.
In analogy to the classical case,
a dg Lie algebra is \emph{complete} if the canonical map to its completion is an isomorphism.
Complete dg Lie algebras form a category $\mathsf{cDGL}$ when 
endowed with chain maps that respect the bracket and the filtration.

\subsubsection{The Lawrence--Sullivan interval, its subdivision, and cDGL models}
\label{sec: LS interval}

The LS (\emph{Lawrence--Sullivan}) interval \cite{LawSul}
is the complete dg Lie algebra 
$$\mathcal L_1 = \left(\widehat{\mathbb{L}}\left(a,b,x\right),\partial\right),$$ 
where $a,b$ are Maurer--Cartan elements and $x$ is a gauge from $b$ to $a$.
That is, $a,b$ are degree $-1$ elements satisfying the Maurer--Cartan equation
$$\partial z + \frac{1}{2}[z,z] = 0,$$ 
and $x$ is a degree $0$ element whose differential codifies a vector field that travels from $b$ to $a$ in unit time,
$$\partial x = \ad_x (b) + \sum_{n\geq 0} \frac{B_n}{n!} \ad_x^n\left(b-a\right).$$ 
Complete dg Lie algebras such as $\mathcal L_1$ can be geometrically realized into topological spaces or simplicial sets by using the realization functor
$$\langle - \rangle : \mathsf{cDGL} \to \mathsf{sSet}$$ 
described in \cite{LieModels}. 
In particular, 
the realization of the Lawrence--Sullivan interval $\mathcal L_1$ 
has the homotopy type of the unit interval together with a disjoint base point \cite{BM1}, 
$$\left\langle \mathcal L_1 \right\rangle \simeq I^+,$$
justifying the name of  $\mathcal L_1$.
The geometric realization of any complete dg Lie algebra as above always carries an artificial external base point that corresponds to the zero Maurer--Cartan element. 
From now on, we drop this external base point from the notation.
There is 
a left adjoint $\mathcal L$ to the realization functor.
Together, these functors form a Quillen pair
\[
\mathcal{L} :  \mathsf{sSet} \leftrightarrows \mathsf{cDGL} :\langle - \rangle
\]
once the classical 
model structure on simplicial sets is transferred to complete dg Lie algebras,
see \cite{LieModels}.

\medskip

Lawrence and Sullivan  
furthermore gave a complete dg Lie model for the first subdivision $\operatorname{sd}_1\left(I\right)$ of the interval \cite[Theorem 2]{LawSul}.
To achieve this,
consider two disjoint LS-intervals 
$$
\left(\widehat{\mathbb{L}}\left(a,b,x\right),\partial\right) \qquad \textrm{and} \qquad \left(\widehat{\mathbb{L}}\left(c,d,y\right),\partial\right).
$$
Take their coproduct in $\mathsf{cDGL}$, 
and take the quotient that identifies the end point $b$ of the
first interval with the end point $c$ of the second one.
Renaming $a$ as $a'$ and $d$ as $b'$, this gives the model of two
glued intervals,
$$\left(\widehat{\mathbb{L}}\left(a',c, b',x, y\right),\partial \right).$$
Here, $a'$, $c$ and $b'$ are Maurer--Cartan elements,
and $x$, $y$ are degree $0$ elements whose differentials are given by
$$\partial x=\ad_x(c)+\sum_{i\geq 0}\frac{B_i}{i!}\ad_x^i\left(c-a'\right)\qquad \textrm{ and }\qquad \partial y=\ad_y\left(b'\right)+\sum_{i\geq 0}\frac{B_i}{i!}\ad_y^i\left(b'-c\right).$$
The orientation is taken as in the following picture.

\vskip .8cm
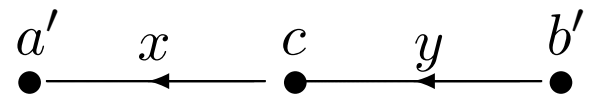
\begin{picture}(0,0)(-130,0)
	\put(0,0){$\bullet$}\put(0,8){$a'$} \put(50, 0){$\bullet$} \put(50,8){$c$}\put(100,0){$\bullet$}\put(100,8){$b'$}
	\put(23,7){$x$}\put(75,7){$y$}
	\qbezier(6,3)(25,3)(47,3)
	\put(45,3){\vector(-1,0){20}}
	\qbezier(52,3)(71,3)(99,3)
	\put(95,3){\vector(-1,0){20}}
\end{picture}
\vskip .5cm

Being more precise,
there is the following result, 
where $\left(\widehat{\mathbb{L}}(a, b,z), \partial \right)$ is a Lawrence--Sullivan interval.

\begin{theorem}
\label{subdivision}
\emph{\cite[Theorem 2]{LawSul}}
The assignment $p(a)=a'$, $p(b)=b'$ and $p(z)=x* y$ defines a
morphism of complete differential graded Lie algebras
$$p\colon \widehat{\mathbb{L}}\left(a,b,z\right)\longrightarrow
\widehat{\mathbb{L}}\left(a',c, b',x, y\right).$$
\end{theorem}

\medskip

We will briefly use in Section \ref{sec: Choice in Bernoulli} the complete dg Lie algebra model of the zero simplex, 
$\mathcal L_0$.
This is the complete dg Lie algebra generated by a single non-trivial Maurer--Cartan element,
$$\mathcal L_0 = \left(\CompLie\left(u\right),\partial\right),$$
where $|u|=-1$ and $\partial u = -\frac{1}{2}[u,u].$
The geometric realization of $\mathcal L_0$ is a discrete two-point space.

\section{The generalized Bernoulli numbers and the $2$-disk}
\label{sec: Number theory}

The gauge action of $\mathcal L_1$, 
the existence of a particular cylinder on the homotopy category of differential graded associative algebras,
and a family of Euler-type identities for Bernoulli numbers are equivalent formulations of the same statement \cite{BCM}, namely the equation
$$\partial^2 = 0.$$
The same equation in the complete dg Lie algebra modelled on the
$2$-disk necessarily yields arithmetical identities, and the question
is what form these identities take.
The Euler-type identities of \cite{BCM} are convolution equations of the form
\begin{equation*} 
\sum_{k=0}^n \lambda_k B_kB_{n-k}=0,
\end{equation*}
where $\lambda_k$ are rational numbers and $B_k$ are the Bernoulli numbers.
Hence, we must first answer the following question:
what should the $2$-dimensional generalization of a Bernoulli number be?
In this section, 
we propose an answer to this question,
and we show that the equation $\partial^2=0$ is equivalent to an
identity involving these new numbers and the classical Bernoulli
numbers, from which we derive explicit quadratic relations among them.

The organization of this section is the following.
In Section \ref{sec: Model of the 2-disk} we give
an explicit complete dg Lie algebra built on a cell decomposition
of the unit $2$-disk,
and describe its universal enveloping algebra.
We define the generalized Bernoulli numbers in 
Section \ref{sec: generalized Bernoulli and 2-disk}.
The idea was already explained in the introduction: 
by analogy with the classical Bernoulli numbers,
these generalized Bernoulli numbers will appear as coefficients in certain expansion formulas for the BCH product in the 
universal enveloping algebra of Section \ref{sec: Model of the 2-disk}.
The identities among the generalized 
Bernoulli numbers are collected in 
Theorem \ref{teo: mainintro}, their equivalence with the equation
$\partial^2=0$ is explained in Remark \ref{rem:equivalence},
and some specific examples are 
computed in Section \ref{sec: examples}.

\subsection{A complete dg Lie algebra for the $2$-disk and its
universal enveloping algebra}
\label{sec: Model of the 2-disk}

Based on the ideas of Section \ref{sec: LS interval},
we now construct a complete dg Lie algebra modelled on the cell
decomposition of the closed $2$-disk with two vertices, two edges
and one $2$-cell. To do so, we glue the end points of the
subdivision $\operatorname{sd}_1\left(I\right)$, obtaining a
circle, and then we attach a $2$-cell along the loop $x*y$. We shall
only use that the resulting derivation squares to zero
(Theorem \ref{disk}). In particular, we do not need, and do not
prove here, that its realization is homotopy equivalent to the
disk. This complete dg Lie algebra appears in \cite[Section 2]{GGL},
where it is obtained from the model of the subdivided interval as in
the proof below. The analogous models of the triangle can be found
in \cite{BFMT4,GL}. We include the proof for completeness.

\begin{theorem}
	\label{disk}
	Let $\widehat{\mathbb{L}}\left(a,b,x,y,\varepsilon\right)$ 
	be the complete free graded Lie algebra 
	generated by $a,b$ in degree $-1$, 
	by $x,y$ in degree $0$,
	and by $\varepsilon$ in degree $1$. 
	Then, the degree $-1$ derivation determined by
	$$\partial a=-\tfrac{1}{2}[a,a],\qquad \partial b=-\tfrac{1}{2}[b,b],\qquad \partial \varepsilon=x* y-\ad_\varepsilon(a),$$
	$$\partial x=\ad_x(b)+\sum_{n\geq 0}\frac{B_n}{n!}\ad_x^n(b-a),\qquad \partial y=\ad_y(a)+\sum_{n\geq 0}\frac{B_n}{n!}\ad_y^n(a-b),$$
	is a differential whose linear term on 
	$\varepsilon$ is 
	$\partial_1 \varepsilon=x+y$.
\end{theorem}

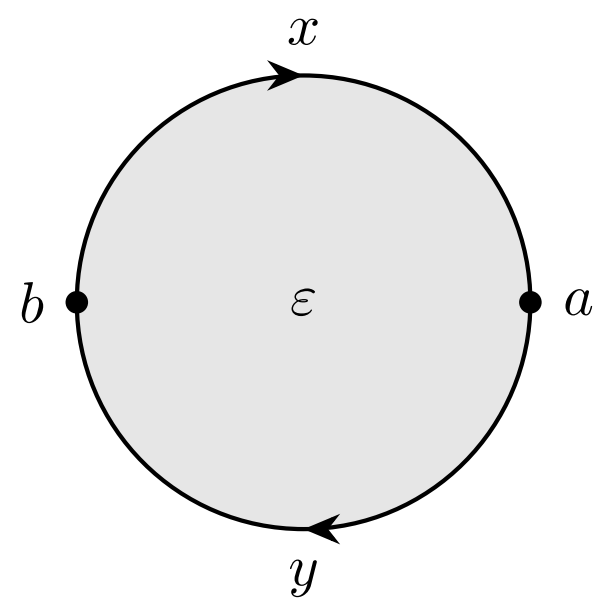
\begin{center}
\begin{tikzpicture}[scale=0.8]
	\pic {celldisk};
	\node at (0,0) {$\varepsilon$};
\end{tikzpicture}
\end{center}

For $z$ of degree $-1$ we write $\partial_z=\partial+\ad_z$, which is
a differential when $z$ is a Maurer--Cartan element. Since $|a|=-1$
and $|\varepsilon|=1$, we have $[a,\varepsilon]=[\varepsilon,a]$, so
the formula for $\partial\varepsilon$ in Theorem \ref{disk} is
equivalent to
$$\partial_a\varepsilon=x*y.$$
That is, the $2$-cell $\varepsilon$ is attached along the loop $x*y$
based at $a$.

\begin{proof}[Proof of Theorem \ref{disk}]
Consider the model of two glued intervals,
$(\widehat{\mathbb{L}}(a',c,b',x,y),\partial)$, and let
$(\widehat{\mathbb{L}}(a,b,z),\partial)$ be a Lawrence--Sullivan
interval. By Theorem \ref{subdivision}, the map $p$ given by $p(a)=a'$, $p(b)=b'$ and $p(z)=x*y$ commutes with
differentials, that is,
$$\partial(x*y)=p(\partial z)=\ad_{x*y}(b')+\sum_{i\geq 0}\frac{B_i}{i!}\ad_{x*y}^i(b'-a').$$
Passing to the quotients by the dg ideals generated by $b-a$ and by
$b'-a'$, and renaming $a'$ as $a$ and $c$ as $b$, the differential of
$\widehat{\mathbb{L}}(a',c,x,y)$ becomes the one of the statement on
$a,b,x,y$, and $\partial(x*y)=\ad_{x*y}(a)$. Then, by the Jacobi
identity,
\[
\partial^2\varepsilon=[x*y,a]-[x*y-[\varepsilon,a],a]-\tfrac{1}{2}[\varepsilon,[a,a]]
=[[\varepsilon,a],a]-\tfrac{1}{2}[\varepsilon,[a,a]]=0. \qedhere
\]
\end{proof}

In the next section,
we work with the universal enveloping algebra of the complete
dg Lie algebra of Theorem \ref{disk}.
It is $\left(\widehat{T}\left(a, b, x, y, \varepsilon\right),
D\right)$, with
$$Da=-a^2, \qquad Db=-b^2, \qquad
D\varepsilon=x* y-\varepsilon a-a\varepsilon,$$
\begin{equation}
\label{1celdas}
Dx=xb-bx+\sum_{p,q\geq 0}\beta_{(p,q)}\,x^p(b-a)x^q,\qquad
Dy=ya-ay+\sum_{p,q\geq 0}\beta_{(p,q)}\,y^p(a-b)y^q,
\end{equation}
where $\beta_{(p,q)}=(-1)^q\frac{B_{p+q}}{p!\,q!}$.
This follows from Theorem \ref{disk} and the identity
$$\ad_z^n(m)=\sum_{p+q=n}\binom{n}{p}(-1)^q z^p m z^q,$$
valid for $z$ of degree $0$.

\subsection{Generalized Bernoulli numbers
	and their identities 
	inherited from the $2$-disk}
\label{sec: generalized Bernoulli and 2-disk}

In this section, 
we work in the complete tensor algebra 
$\widehat T\left(x,y\right)$ on two degree zero generators,
making heavy use of the BCH formula for  $x*y$ explained in Section \ref{sec: BCH prerequisites}. 
Recall that
$$x* y = \sum_{n\geq 1} z_n(x,y),$$
where 
$$z_n(x,y) = \sum_{k=1}^n\frac{(-1)^{k-1}}{k}\sum \frac{x^{p_1}y^{q_1}\cdots x^{p_k}y^{q_k}}{p_1!q_1!\cdots p_k!q_k!}.$$
The second sum above ranges over all possible partitions of $n$ by non-negative integers $$(p_1,...,p_k,q_1,...,q_k)$$ 
satisfying $p_i + q_i \geq 1$ for all $i$, 
so that $p_1+q_1+\cdots +p_k+q_k=n.$

Now, let $w$ be a word in $x$ and $y$.
That is, 
$$w = x^{h_1}y^{k_1}\cdots x^{h_n}y^{k_n}, \quad \textrm{$n\geq 1$,}$$
where
\begin{itemize}
	\item $h_1 \geq 0$ and $h_i\geq 1$ for $i=2,\dots ,n$,
	\item $k_i\geq 1$ for $i=1,\dots, n-1$ and $k_n\geq 0$.
\end{itemize}  
By convention,  $x^0 = y^0 = 1$ are omitted from the word, which can therefore begin and/or end with $x$ and/or $y$.
For convenience, we identify the word $w$ with the vector of integers 
\begin{equation}
	\label{ecu: vector of integers}
	\left( \left(h_1, k_1\right),\dots , \left(h_n, k_n\right)\right).
\end{equation}
Let $c_w=c_{\left( \left(h_1, k_1\right),\dots , \left(h_n, k_n\right)\right)} \in \mathbb{Q}$ be the coefficient of the word $w$ in the expression of $x* y$. 
The \emph{generalized Bernoulli number}\footnote{Not to be confused
with the generalized Bernoulli numbers $B_{n,\chi}$ attached to
Dirichlet characters, see for instance \cite{BerAndZeta}.}
is the rational number
\begin{equation*}
\mathcal{B}_w=h_1!k_1!\cdots h_n!k_n!c_{\left( \left(h_1, k_1\right),\dots , \left(h_n, k_n\right)\right)}.
\end{equation*}
This rational number generalizes the (signed) Bernoulli numbers.
Indeed, for all $p,q\geq 0$ with
$p+q\geq 1$ we have
$$\mathcal{B}_{x^pyx^q}=\left(-1\right)^p B_{p+q},$$
because the part of $x*y$ linear in $y$ is
$\sum_{n\geq 0}(-1)^n\frac{B_n}{n!}\ad_x^n(y)$, see for instance
\cite{BonfiglioliFulci}.

Moreover, 
let $w_1=x^{r_1}y^{s_1}\cdots x^{r_{n}}y^{s_{n}}$ and $w_2=x^{t_1}y^{u_1}\cdots x^{t_{m}}y^{u_{m}}$
be two words on $x$ and $y$. We introduce two different operations to be performed with words.
\begin{itemize}
	\item The first one is the
	{\em product} $w_1w_2$,
	denoted by juxtaposition.
	This is simply the (associative, unital) multiplication on the free monoid generated by $x$ and $y$.
	More precisely, 
	and in order to set up the notation, 
	the product $w_1w_2$ is
	\begin{equation}
		\label{yuxtaposicion}
		x^{r_1}y^{s_1}\cdots x^{r_{n}}y^{s_{n}}x^{t_1}y^{u_1}\cdots x^{t_{m}}y^{u_{m}}
	\end{equation}
	if $s_n\not=0$ and $t_1\not=0$. If on the contrary $s_n=t_1=0$,
	then the terms $y^{s_n}$ and $x^{t_1}$ are dropped. 
	But if one, and only one,
	of the two exponents $s_n$ and $t_1$ is zero, 
	then the word obtained from the expression 
	$(\ref{yuxtaposicion})$ 
	(which in that case is not a correct word) 
	is obtained by replacing $x^{r_n}y^{s_n}x^{t_1}$ with $x^{r_n+t_1}$ for $s_n=0$, 
	or  $y^{s_n}x^{t_1}y^{u_1}$ by $y^{s_n+u_1}$ for $t_1=0$.

	\item The second operation  is defined only when the two words $w_1,w_2$ have the same number of pairs,
	i.e., when $n=m$.
	We denote it $w_1+w_2$ and call it  the \emph{sum}.
	If 
	$w_1=x^{r_1}y^{s_1}\cdots x^{r_{n}}y^{s_{n}}$ and $w_2=x^{t_1}y^{u_1}\cdots x^{t_{n}}y^{u_{n}}$, then
	$$w_1+w_2=x^{r_1+t_1}y^{s_1+u_1}\cdots x^{r_{n}+t_n}y^{s_{n}+u_n}.$$
	This operation is simply the addition of the vectors that correspond to $w_1$ and $w_2$,
	recall (\ref{ecu: vector of integers}).
In this operation,
	we can eventually consider vectors with more $0$'s than in a regular word $w$ where only $h_1$ and $k_n$ are allowed to be $0$. Then the notation $\underbrace{(h, k)}_i$ stands for 
	$$\left( (0,0),\dots ,\underbrace{(h,k)}_i, \dots ,(0,0)\right).$$ 
	The number of zeroes after the $i$th coordinate is irrelevant.
\end{itemize} 
We will use  the following notation for multinomial numbers: 
$\binom{n}{ m_1, \ m_2}=\frac{n!}{m_1!m_2!(n-m_1-m_2)!}$.
Our main result reads as follows.

\begin{theorem}
	\label{teo: mainintro}
For any pair 
$(w_1, w_2)=\left( \left(x^{r_1}y^{s_1}\cdots x^{r_{n}}y^{s_{n}}\right),\left(x^{t_1}y^{u_1}\cdots x^{t_{m}}y^{u_{m}}\right)\right)$ of words on $x$ and $y$, the following identities hold:
\begin{enumerate}
\item If $s_n\neq 0$ and $t_1\neq 0$, then
\begin{align*}
&\sum_{\ell=0}^{t_1}(-1)^{\ell}\binom{s_n+t_1+1}{\ell,\ s_n}
B_\ell\mathcal{B}_{w_1w_2+\underbrace{(1-\ell, 0)}_{n+1}}\\
&\qquad-\sum_{\ell=0}^{s_n}(-1)^\ell\binom{s_n+t_1+1}{\ell,\ t_1}
B_\ell\mathcal{B}_{w_1w_2+\underbrace{(0,1-\ell)}_{n}}=0.
\end{align*}

\item If $s_n=t_1=0$, so that
$w_1w_2=x^{r_1}y^{s_1}\cdots x^{r_n}y^{u_1}x^{t_2}\cdots
x^{t_m}y^{u_m}$ has $n+m-1$ pairs, then
\begin{align*}
&\sum_{\ell=0}^{r_n}(-1)^{\ell}\binom{r_n+u_1+1}{\ell,\ u_1}
B_\ell\mathcal{B}_{w_1w_2+\underbrace{(1-\ell,0)}_{n}}\\
&\qquad-\sum_{\ell=0}^{u_1}(-1)^\ell\binom{r_n+u_1+1}{\ell,\ r_n}
B_\ell\mathcal{B}_{w_1w_2+\underbrace{(0,1-\ell)}_{n}}=0.
\end{align*}
\item Let $z\in\{x,y\}$ and let $\bar z$ be the other letter. If
$w_1=w_1'z^{\mu}$ and $w_2=z^{\nu}w_2'$ with $\mu,\nu\geq 1$, where
$w_1'$ does not end and $w_2'$ does not begin with $z$, then, with
$N=\mu+\nu+1$,
\begin{align*}
\mathcal{B}_{w_1'z^{N}w_2'}
+\sum_{\substack{2\leq\ell\leq \mu+\nu\\ \ell\ \mathrm{even}}}
\binom{N}{\ell}\lambda_\ell(\mu,\nu)\,B_\ell\,
\mathcal{B}_{w_1'z^{N-\ell}w_2'}
=\frac{N!}{\mu!\,\nu!}\,\mathcal{B}_{w_1\bar zw_2},
\end{align*}
where $\lambda_\ell(\mu,\nu)=[\ell>\mu](-1)^{\mu}\binom{\ell-1}{\mu}
+[\ell>\nu](-1)^{\nu}\binom{\ell-1}{\nu}$ and $[P]$ is $1$ if $P$
holds and $0$ otherwise.
\end{enumerate}
\end{theorem}

\begin{proof}
Write $\Psi=x*y=\sum_w c_w\,w$, the sum running over all
non-empty words $w$ in $x$ and $y$, and for $p,q\geq 0$ put
$$\alpha_{(p,q)}=(-1)^{p+1}\frac{B_{p+q}}{p!\,q!},\qquad
\beta_{(p,q)}=(-1)^{q}\frac{B_{p+q}}{p!\,q!}.$$

\emph{Step 1.} We have
\begin{equation}\label{especialy}
Dx=xa-ax+\Phi_x,\qquad Dy=ya-ay+\Phi_y,
\end{equation}
where $\Phi_x=\sum_{p,q\geq 0}\alpha_{(p,q)}x^p(a-b)x^q$ and
$\Phi_y=\sum_{p,q\geq 0}\beta_{(p,q)}y^p(a-b)y^q$.
The second identity is (\ref{1celdas}). For the first,
substitute $xb-bx=xa-ax-x(a-b)+(a-b)x$ in (\ref{1celdas}). Then
the coefficient of $x^p(a-b)x^q$ becomes $-\beta_{(p,q)}$, corrected
by $-1$ for $(p,q)=(1,0)$ and by $+1$ for $(p,q)=(0,1)$, and this
equals $\alpha_{(p,q)}$ because $B_1=-\frac{1}{2}$ and
$B_{2k+1}=0$ for $k\geq 1$.

\emph{Step 2.} Since $D$ is a derivation of degree $-1$ and
$Da=-a^2$,
$$D^2\varepsilon=D\Psi-(D\varepsilon)a+\varepsilon Da-(Da)\varepsilon
+aD\varepsilon=D\Psi-\Psi a+a\Psi.$$
As $D^2=0$ by Theorem \ref{disk},
\begin{equation}\label{equivalente}
D\Psi=\Psi a-a\Psi.
\end{equation}

\emph{Step 3.} By (\ref{especialy}), for $z\in\{x,y\}$ and
$m\geq 1$ we have
$$D(z^m)=z^ma-az^m+\sum_{i+j=m-1}z^i\Phi_zz^j.$$
Hence, for every
word $w$, $D(w)=wa-aw+\Gamma_w$, where $\Gamma_w$ is the sum of the
elements obtained from $w$ by replacing one occurrence of a letter
$z\in\{x,y\}$ by $\Phi_z$. By (\ref{equivalente}),
\begin{equation}\label{laecuacion}
\sum_w c_w\,\Gamma_w=0.
\end{equation}
Every summand of $\Gamma_w$ is a multiple of an element
$u(a-b)v$ with $u,v$ words in $x$ and $y$, and these elements are
linearly independent. Hence, for every pair $(w_1,w_2)$ of
non-empty words, the coefficient of $w_1(a-b)w_2$ in the left-hand
side of (\ref{laecuacion}) vanishes.

\emph{Step 4.} Fix non-empty words $w_1,w_2$ and write
$w_1=w_1'z_1^{\,\mu}$ and $w_2=z_2^{\,\nu}w_2'$ with
$z_1,z_2\in\{x,y\}$ and $\mu,\nu\geq 1$, where $w_1'$ is empty or ends
with a letter different from $z_1$, and $w_2'$ is empty or begins
with a letter different from $z_2$. A summand of $\Gamma_w$ is a
multiple of $w_1(a-b)w_2$ exactly when $w=uzv$ with $z\in\{x,y\}$
and the summand is $\alpha_{(p,q)}\,uz^p(a-b)z^qv$ (if $z=x$) or
$\beta_{(p,q)}\,uz^p(a-b)z^qv$ (if $z=y$), with $uz^p=w_1$ and
$z^qv=w_2$. We list these possibilities.

$(i)$ $z_1=y$, $z_2=x$. Either $z=x$, $p=0$ and
$w=w_1x^{\nu+1-q}w_2'$ with $0\leq q\leq \nu$, or $z=y$, $q=0$ and
$w=w_1'y^{\mu+1-p}w_2$ with $0\leq p\leq \mu$. Since
$\alpha_{(0,\ell)}=-B_\ell/\ell!$ and
$\beta_{(\ell,0)}=B_\ell/\ell!$,
\begin{equation}\label{eq:caseA}
\sum_{\ell=0}^{\nu}\frac{B_\ell}{\ell!}\,c_{w_1x^{\nu+1-\ell}w_2'}
=\sum_{\ell=0}^{\mu}\frac{B_\ell}{\ell!}\,c_{w_1'y^{\mu+1-\ell}w_2}.
\end{equation}

$(ii)$ $z_1=x$, $z_2=y$. Either $z=x$, $q=0$ and
$w=w_1'x^{\mu+1-p}w_2$, or $z=y$, $p=0$ and $w=w_1y^{\nu+1-q}w_2'$.
Since $\alpha_{(\ell,0)}=(-1)^{\ell+1}B_\ell/\ell!$ and
$\beta_{(0,\ell)}=(-1)^{\ell}B_\ell/\ell!$,
\begin{equation}\label{eq:caseB}
\sum_{\ell=0}^{\mu}(-1)^\ell\frac{B_\ell}{\ell!}\,
c_{w_1'x^{\mu+1-\ell}w_2}
=\sum_{\ell=0}^{\nu}(-1)^\ell\frac{B_\ell}{\ell!}\,
c_{w_1y^{\nu+1-\ell}w_2'}.
\end{equation}

$(iii)$ $z_1=z_2=x$. Either $z=x$, $0\leq p\leq \mu$,
$0\leq q\leq \nu$ and $w=w_1'x^{\mu+\nu+1-p-q}w_2'$, or $z=y$, $p=q=0$
and $w=w_1yw_2$. Since $\beta_{(0,0)}=1$,
\begin{equation}\label{eq:caseC}
\sum_{p=0}^{\mu}\sum_{q=0}^{\nu}(-1)^{p}\frac{B_{p+q}}{p!\,q!}\,
c_{w_1'x^{\mu+\nu+1-p-q}w_2'}=c_{w_1yw_2}.
\end{equation}

$(iv)$ $z_1=z_2=y$. In the same way, since $\alpha_{(0,0)}=-1$,
\begin{equation}\label{eq:caseD}
\sum_{p=0}^{\mu}\sum_{q=0}^{\nu}(-1)^{q}\frac{B_{p+q}}{p!\,q!}\,
c_{w_1'y^{\mu+\nu+1-p-q}w_2'}=c_{w_1xw_2}.
\end{equation}

\emph{Step 5.} For a word $w$ let $P(w)$ be the product of the
factorials of the exponents of $w$, so that
$c_w=\mathcal{B}_w/P(w)$, and put $N=\mu+\nu+1$.

In case $(i)$ we have $\mu=s_n$, $\nu=t_1$,
$w_1x^{\nu+1-\ell}w_2'=w_1w_2+\underbrace{(1-\ell,0)}_{n+1}$ and
$w_1'y^{\mu+1-\ell}w_2=w_1w_2+\underbrace{(0,1-\ell)}_{n}$.
Multiplying (\ref{eq:caseA}) by $P(w_1w_2)N!/(\mu!\,\nu!)$ gives
$$\sum_{\ell=0}^{t_1}\binom{N}{\ell,\ s_n}B_\ell\,
\mathcal{B}_{w_1w_2+\underbrace{(1-\ell,0)}_{n+1}}
-\sum_{\ell=0}^{s_n}\binom{N}{\ell,\ t_1}B_\ell\,
\mathcal{B}_{w_1w_2+\underbrace{(0,1-\ell)}_{n}}=0,$$
which is identity (1) without the signs $(-1)^\ell$. The two forms
are equivalent, since $B_\ell=0$ for odd $\ell\geq 3$ and the two
terms with $\ell=1$, each equal to
$\frac{N!}{s_n!\,t_1!}B_1\mathcal{B}_{w_1w_2}$, cancel in either
form.

In case $(ii)$ we have $\mu=r_n$, $\nu=u_1$ and
$w_1w_2=w_1'x^{r_n}y^{u_1}w_2'$, so that
$w_1'x^{\mu+1-\ell}w_2=w_1w_2+\underbrace{(1-\ell,0)}_{n}$ and
$w_1y^{\nu+1-\ell}w_2'=w_1w_2+\underbrace{(0,1-\ell)}_{n}$.
Multiplying (\ref{eq:caseB}) by the factor
$P(w_1w_2)N!/(\mu!\,\nu!)$ gives identity (2).

In cases $(iii)$ and $(iv)$ let $z=z_1=z_2$ and let $\bar z$ be the
other letter. For $0\leq\ell\leq \mu+\nu$ put
$S_\ell=\sum(-1)^p\binom{\ell}{p}$, the sum over the $p$ with
$\max(0,\ell-\nu)\leq p\leq\min(\mu,\ell)$. From
$\sum_{p=0}^{b}(-1)^p\binom{\ell}{p}=(-1)^b\binom{\ell-1}{b}$,
valid for $0\leq b<\ell$, one gets $S_0=1$, $S_1=0$, and, for even
$\ell\geq 2$,
$$S_\ell=\lambda_\ell(\mu,\nu):=[\ell>\mu]\,(-1)^\mu\binom{\ell-1}{\mu}
+[\ell>\nu]\,(-1)^\nu\binom{\ell-1}{\nu},$$
where $[P]$ is $1$ if $P$ holds and $0$ otherwise. For even $\ell$
and for $\ell=1$, the same values are obtained with $(-1)^q$ in place
of $(-1)^p$. Multiplying (\ref{eq:caseC}), or (\ref{eq:caseD}), by
$P(w_1')P(w_2')N!$ and using $B_\ell=0$ for odd $\ell\geq 3$, we
obtain
\begin{equation}\label{eq:mixed}
\mathcal{B}_{w_1'z^{N}w_2'}
+\sum_{\substack{2\leq\ell\leq \mu+\nu\\ \ell\ \mathrm{even}}}
\binom{N}{\ell}\lambda_\ell(\mu,\nu)\,B_\ell\,
\mathcal{B}_{w_1'z^{N-\ell}w_2'}
=\frac{N!}{\mu!\,\nu!}\,\mathcal{B}_{w_1\bar zw_2}.
\end{equation}
This completes the proof.
\end{proof}
\begin{remark}\label{rem:equivalence}
Identity (\ref{laecuacion}) involves the generalized Bernoulli
numbers $\mathcal{B}_w=P(w)c_w$ and, through $\Phi_x$ and $\Phi_y$,
the classical ones. It is equivalent to the existence of the
differential of Theorem \ref{disk}, that is, to the equation
$\partial^2=0$. Indeed, $D^2=\frac{1}{2}[D,D]$ is a derivation which
vanishes on $a$, $b$, $x$ and $y$ (see the proof of
Theorem \ref{disk}), and Steps 2 and 3 of the proof above, which do
not use that $D^2=0$, give
$$D^2\varepsilon=D\Psi-\Psi a+a\Psi=\sum_w c_w\,\Gamma_w.$$
Hence $D^2=0$ if and only if (\ref{laecuacion}) holds. The identities
of Theorem \ref{teo: mainintro} are the components of
(\ref{laecuacion}) on the elements $w_1(a-b)w_2$ with $w_1$ and $w_2$
non-empty, written in terms of the numbers $\mathcal{B}_w$.
\end{remark}
\begin{remark}\label{rem:associativity}
Identity (\ref{laecuacion}), and with it
Theorem \ref{teo: mainintro}, can also be obtained without the
model of the disk. Let $u$ be a further variable of degree $0$, and
work in $\widehat T(x,y,u)$ modulo the closed ideal of elements which are
at least quadratic in $u$. With $g(s)=s/(e^s-1)$, the parts of
$x*u$ and of $(-u)*y$ linear in $u$ are $g(-\ad_x)(u)$ and
$-g(\ad_y)(u)$. The first is the part of $x*y$ linear in $y$,
computed in Section \ref{sec: generalized Bernoulli and 2-disk},
with $u$ in place of $y$. The second follows from the first applied
to $(-y)*u$, because $(-u)*y=-((-y)*u)$. This last identity holds
since $(e^{-u}e^{y})^{-1}=e^{-y}e^{u}$. As $*$ is associative and
$u*(-u)=0$, we have
$(x*u)*((-u)*y)=x*y$, and expanding the left-hand side gives
$\sum_wc_w\widetilde\Gamma_w=0$, where $\widetilde\Gamma_w$ is
obtained from $w$ by replacing, one at a time, each letter $x$ by
$g(-\ad_x)(u)$ and each letter $y$ by $-g(\ad_y)(u)$. Under
$u\mapsto b-a$ these become $\Phi_x$ and $\Phi_y$, and
$\widetilde\Gamma_w$ becomes $\Gamma_w$.
\end{remark}

\subsubsection{Examples}
\label{sec: examples}
In this subsection, we 
discuss some examples and relevant cases of the family of identities of Theorem  \ref{teo: mainintro}.

\begin{example}
	Consider the pair of words 
	$\left( \left(xy^2\right), (xy)\right)$. 
Then, $r_1=1, s_1=2$ and $t_1=1, u_1=1$. 
The corresponding identity in Theorem \ref{teo: mainintro} becomes
\begin{eqnarray*}
&\sum_{\ell =0}^1 (-1)^\ell \binom{4}{\ell,\ 2} B_\ell \mathcal{B}_{xy^2xy+\underbrace{(1-\ell, 0)}_2}&\\
-&\sum_{\ell =0}^2 (-1)^\ell \binom{4}{\ell,\ 1} B_\ell \mathcal{B}_{xy^2xy+\underbrace{(0,1-\ell)}_1}&=0.
\end{eqnarray*}
The left hand side of this expression is
\begin{eqnarray*}
\left( 6B_0\mathcal{B}_{xy^2x^2y}-12B_1 \mathcal{B}_{xy^2xy}\right)-\left( 4 B_0\mathcal{B}_{xy^3xy}-12 B_1\mathcal{B}_{xy^2xy}+12B_2\mathcal{B}_{xyxy}\right)&=\\
\left( 6\mathcal{B}_{xy^2x^2y}+6 \mathcal{B}_{xy^2xy}\right)-\left( 4\mathcal{B}_{xy^3xy}+6\mathcal{B}_{xy^2xy}+2\mathcal{B}_{xyxy}\right)&=\\
6\mathcal{B}_{xy^2x^2y}-4\mathcal{B}_{xy^3xy}-2\mathcal{B}_{xyxy},&
\end{eqnarray*}
since $B_0=1$, $B_1=-\frac{1}{2}$, $B_2=\frac{1}{6}$. 
To compute the remaining terms we use the following closed formula
for $\mathcal{B}_w$. It is a reformulation of Goldberg's theorem
\cite[Theorem 1]{Goldberg}, obtained by writing the polynomials in
Goldberg's integral formula in terms of Stirling numbers of the
second kind and integrating term by term.

If $w=x^{h_1}y^{k_1}\cdots x^{h_n}y^{k_n}$ and  $N=\sum_{j=1}^n \left(h_j+k_j\right)$, 
then
\begin{equation}
	\label{BCHpalabra}
\mathcal{B}_w=\frac{1}{d+1}\sum_{k=d+n}^N(-1)^{k-d-1}\frac{G_k(w)}{\binom{k}{d+1}},
\end{equation}
where
$$d=
\left\{
\begin{array}{ll}
n & \mbox{if } h_1\neq 0 \mbox{ and } k_n\neq 0, \\
n-1 & \mbox{if exactly one of } h_1, k_n \mbox{ is } 0, \\
n-2 & \mbox{if } h_1=k_n=0,
\end{array}
\right.$$

\begin{equation}
G_k(w)=\sum_{\sum_i (\sigma_i+\tau_i)=k}\ \prod_{i=1}^n \sigma_i!\,\tau_i!\,
\stirling{k_i}{\sigma_i}\stirling{h_i}{\tau_i},
\end{equation}
and the symbol $\stirling{a}{b}$ denotes the Stirling number of the
second kind, i.e. the number of ways to partition a set of $a$
objects into $b$ non-empty subsets, with $\stirling{0}{0}=1$.
Then, for example, for the word $w=xy^3xy$, we have $n=2$, $d=2$,
$N=6$, $h_1=1$, $k_1=3$, $h_2=1$, $k_2=1$.
Thus,
\begin{align*}
	\mathcal{B}_{xy^3xy} &= \frac{1}{2+1}\sum_{k=4}^6 (-1)^{k-3}\frac{G_k(xy^3xy)}{\binom{k}{3}}\\ 
	&=\frac{1}{3}\left( -\frac{1}{4}G_4(xy^3xy)+\frac{1}{10}G_5(xy^3xy)-\frac{1}{20}G_6(xy^3xy) \right)=\frac{1}{60},
\end{align*}
since
$$
	G_4(xy^3xy)=\stirling{1}{1}\stirling{3}{1}\stirling{1}{1}\stirling{1}{1}=1,\quad
	G_5(xy^3xy)=2!\stirling{1}{1}\stirling{3}{2}\stirling{1}{1}\stirling{1}{1}=6
$$
and
$$
	G_6(xy^3xy)=3!\stirling{1}{1}\stirling{3}{3}\stirling{1}{1}\stirling{1}{1}=6.
$$
Similar calculations give $\mathcal{B}_{xyxy}=-\frac{1}{12}$ and $\mathcal{B}_{xy^2x^2y}=-\frac{1}{60}$.
Then,
$$6\mathcal{B}_{xy^2x^2y}-4\mathcal{B}_{xy^3xy}-2\mathcal{B}_{xyxy}=-\frac{1}{10}-\frac{1}{15}+\frac{1}{6}=0.$$
\hfill $\square$
\end{example}

\begin{example}\label{ex:euler}
The identities of Theorem \ref{teo: mainintro} contain Euler's
identity
$$(n+1)B_n+\sum_{k=2}^{n-2}\binom{n}{k}B_kB_{n-k}=0,
\qquad n\geq 4 \mbox{ even}.$$
To see this, let $q\geq 3$ be odd, put $n=q+1$, and apply
Theorem \ref{teo: mainintro}(1) to $(w_1,w_2)=(y,x^q)$, that is,
$r_1=0$, $s_1=1$, $t_1=q$, $u_1=0$. Then
$w_1w_2+\underbrace{(1-\ell,0)}_2=yx^{n-\ell}$ and
$w_1w_2+\underbrace{(0,1-\ell)}_1=y^{2-\ell}x^{n-1}$, and the
identity reads
$$\sum_{\ell=0}^{n-1}(-1)^\ell\binom{n+1}{\ell,\ 1}B_\ell\,
\mathcal{B}_{yx^{n-\ell}}
-\sum_{\ell=0}^{1}(-1)^\ell\binom{n+1}{\ell,\ n-1}B_\ell\,
\mathcal{B}_{y^{2-\ell}x^{n-1}}=0.$$
Here $\binom{n+1}{\ell,\ 1}=(n+1)\binom{n}{\ell}$,
$\binom{n+1}{0,\ n-1}=\frac{(n+1)n}{2}$,
$\binom{n+1}{1,\ n-1}=(n+1)n$, and $\mathcal{B}_{yx^m}=B_m$ by the
formula for $\mathcal{B}_{x^pyx^q}$ in
Section \ref{sec: generalized Bernoulli and 2-disk}. Since
$B_{n-1}=0$, all the terms with odd $\ell$ vanish, and dividing by
$n+1$ we obtain
\begin{equation}
\label{ecuacioncilla}
\sum_{\ell=0}^{n-1}\binom{n}{\ell}B_\ell B_{n-\ell}
-\frac{n}{2}\mathcal{B}_{y^2x^{n-1}}=0.
\end{equation}
The terms $\ell=1$ and $\ell=n-1$ of the sum vanish, so the sum is
$B_n+\sum_{\ell=2}^{n-2}\binom{n}{\ell}B_\ell B_{n-\ell}$. On the
other hand, the closed formula for the coefficients of the words
$x^{N-k}y^k$ in $x*y$ given in \cite[Proposition 2]{Hof}, together
with $c_w=(-1)^{N+1}c_{\bar w}$, where $N$ is the length of $w$ and
$\bar w$ is obtained from $w$ by interchanging $x$ and $y$, gives
$\mathcal{B}_{y^2x^{n-1}}=-2B_n-B_{n-1}=-2B_n$. The relation between
$c_w$ and $c_{\bar w}$ holds because $x*y=-((-y)*(-x))$. Hence
$\frac{n}{2}\mathcal{B}_{y^2x^{n-1}}=-nB_n$, and
(\ref{ecuacioncilla}) is Euler's identity.
\hfill $\square$
\end{example}

\section{Symmetries and Bernoulli numbers}
\label{sec: symmetries}

This section is devoted to relating geometric properties of the
objects under study, namely the interval and the $2$-disk, to some
notable facts or peculiarities of the Bernoulli numbers, such as the
existence of two choices for the first Bernoulli number $B_1$, the
vanishing of the odd Bernoulli numbers $B_n$ for $n\geq 3$, or the
equality $B_4=B_8$, which is the only coincidence between two
non-zero Bernoulli numbers with different indices.

The organization of this section is the following.
First, 
we revisit two classical results for the usual Bernoulli
numbers in terms of the topological
properties of the unit interval.
To be precise, 
we show
that the vanishing of the odd Bernoulli numbers except for the
first term $B_1$ is exactly the condition for the symmetry of the
unit interval to be compatible with the Lawrence--Sullivan
differential (Section \ref{sec: vanishing of odd Bernoulli}).
We also prove
that the two distinct choices for the initial term $B_1$ of the sequence 
with opposite sign, namely $\frac{1}{2}$ and $-\frac{1}{2}$, correspond to the two
opposite orientations for the interval
(Section 
\ref{sec: Choice in Bernoulli}).
These are elementary, yet interesting, results.
Then, in Section \ref{sec: disk symmetries}, we consider the
symmetries of the cell decomposition of the $2$-disk of
Section \ref{sec: Model of the 2-disk}: the reflection that
interchanges the two edges, which corresponds to the inversion formula
of the BCH product, and the rotation by $\pi$, which implies the
equality $B_4=B_8$.
\subsection{The vanishing of the odd Bernoulli numbers}
\label{sec: vanishing of odd Bernoulli}

The interval has a symmetry exchanging its end points. On the
Lawrence--Sullivan interval it is given by the automorphism
$\sigma$ of the complete graded Lie algebra $\CompLie(a,b,x)$
determined by $\sigma(a)=b$, $\sigma(b)=a$ and $\sigma(x)=-x$.

\begin{proposition}\label{prop:symmetry}
Let $f(t)=\sum_{n\geq 0}\lambda_nt^n\in\mathbb{Q}[[t]]$ and let
$\partial_f$ be the derivation of degree $-1$ of $\CompLie(a,b,x)$
given by $\partial_fa=-\frac{1}{2}[a,a]$,
$\partial_fb=-\frac{1}{2}[b,b]$ and
$$\partial_fx=\ad_x(b)+\sum_{n\geq 0}\lambda_n\ad_x^n(b-a).$$
Then $\sigma\partial_f=\partial_f\sigma$ if and only if
$f(t)-f(-t)=-t$, that is, if and only if $\lambda_1=-\frac{1}{2}$
and $\lambda_{2k+1}=0$ for all $k\geq 1$.
\end{proposition}

\begin{proof}
Both maps agree on $a$ and $b$. On $x$,
$$\sigma(\partial_fx)=-\ad_x(a)+f(-\ad_x)(a-b),\qquad
\partial_f(\sigma x)=-\partial_fx=-\ad_x(b)+f(\ad_x)(a-b),$$
so that
$\sigma(\partial_fx)-\partial_f(\sigma x)
=\bigl(f(-\ad_x)-f(\ad_x)-\ad_x\bigr)(a-b)$. The elements
$\ad_x^n(a-b)$, $n\geq 0$, are non-zero and have different degrees
in $x$, hence are linearly independent, and the claim follows by
comparing coefficients.
\end{proof}

For the Lawrence--Sullivan interval $\mathcal L_1$ we have
$\lambda_n=B_n/n!$, and the identity
$\frac{t}{e^t-1}-\frac{-t}{e^{-t}-1}=-t$ shows that $\sigma$ is an
automorphism of $\mathcal L_1$, which reverses the gauge $x$, in
accordance with the fact that $-x$ is a gauge from $a$ to $b$
\cite[Remark, p.8]{LawSul}. By Proposition \ref{prop:symmetry},
the conditions $B_1=-\frac{1}{2}$ and $B_{2k+1}=0$ for $k\geq 1$
are exactly what makes this symmetry of the interval compatible
with the Lawrence--Sullivan differential.

\subsection{The two choices of first Bernoulli number}
\label{sec: Choice in Bernoulli}

The realization of the Lawrence--Sullivan interval carries the
external base point given by the Maurer--Cartan element $0$
(Section \ref{sec: LS interval}). Identifying an end point of the
interval with it gives the quotients
$$\CompLie\left(a,b,x\right)\twoheadrightarrow\CompLie\left(a,x\right)
\qquad \textrm{ and } \qquad
\CompLie\left(a,b,x\right)\twoheadrightarrow\CompLie\left(b,x\right)$$
by the ideals generated by $b$ and by $a$, respectively. These ideals
are dg ideals, since $\partial b=-\frac{1}{2}[b,b]$ and
$\partial a=-\frac{1}{2}[a,a]$. The two quotients are the
cofibers in $\mathsf{cDGL}$ of the morphisms
$i_b,i_a\colon\mathcal L_0=\CompLie(u)\to\CompLie(a,b,x)$ given by
$i_b(u)=b$ and $i_a(u)=a$. We regard $\CompLie(a,x)$ and
$\CompLie(b,x)$ as models of the interval based at $b$ and at $a$,
respectively.
Let us look at the induced differentials on the respective cofibers.

\medskip

\noindent $(i)$ 
On the algebra
$\CompLie\left(a,x\right)$, 
since $b$ was identified with the zero Maurer--Cartan element, and the gauge $x$ was from $b$ to $a$, 
we get 
$$\partial x = \sum_{i\geq 0} \frac{B_i}{i!}\ad_x^i(-a).$$

\noindent $(ii)$
On the algebra
$\CompLie\left(b,x\right)$, 
it was $a$ that was identified with the zero Maurer--Cartan element. 
Thus, 
\begin{align*}
	\partial x &= \ad_x(b) + \sum_{i\geq 0}\frac{B_i}{i!}\ad_x^i(b)
	= [x,b] + b - \frac{1}{2} [x,b] + \sum_{i\geq 2} \frac{B_i}{i!}\ad_x^i(b) \\
	&= b+\frac{1}{2} [x,b] + \sum_{i\geq 2} \frac{B_i}{i!}\ad_x^i(b).
\end{align*}
It is possible to express the differential $\partial x$ above in a more compact form.
To do so, 
it suffices to change the first Bernoulli
number from $B_1 = -\frac{1}{2}$ to $B_1 =\frac{1}{2}$.
With this change in mind, 
we can write the differential on $(ii)$ as
$$\partial x = \sum_{i\geq 0} \frac{B_i}{i!}\ad_x^i(b).$$

Summarizing, let $u$ be a Maurer--Cartan element and $x$ a
generator of degree $0$, and write $B^{+}_i$ for the Bernoulli
numbers with the choice $B^{+}_1=\frac{1}{2}$ and $B^{-}_i=B_i$
for those with $B^{-}_1=-\frac{1}{2}$, so that
$B^{+}_i=(-1)^iB^{-}_i$ for all $i\geq 0$. By $(ii)$, the complete
dg Lie algebra $\CompLie\left(u,x\right)$ with
$$\partial x = \sum_{i\geq 0} \frac{B^{+}_i}{i!}\ad_x^i(u)$$
represents the interval oriented from $u$ to $0$, since $x$ is a
gauge from $u$ to $0$. By $(i)$, with $u=a$, the one with
$$\partial x = -\sum_{i\geq 0} \frac{B^{-}_i}{i!}\ad_x^i(u)$$
represents the interval oriented from $0$ to $u$. The substitution
$x\mapsto -x$, which reverses the orientation, interchanges the two
differentials. Thus the two choices of the first Bernoulli number
correspond to the two orientations of the interval. The global sign
in the second formula cannot be dropped. Indeed, for
$\partial x=\sum_{i\geq 0}\frac{B^{-}_i}{i!}\ad_x^i(u)$ the
component of $\partial^2x$ in $\CompLie(u)$ is $-[u,u]\neq 0$.

\subsection{Symmetries of the \texorpdfstring{$2$}{2}-disk}
\label{sec: disk symmetries}

\begin{figure}[ht]
\centering
\begin{tikzpicture}
	\begin{scope}
		\pic {celldisk};
		\node at (-0.5,0) {$\varepsilon$};
		\draw[symmetry, {Stealth[length=2.6mm]}-{Stealth[length=2.6mm]}]
			(0,-0.95) -- (0,0.95);
		\node at (0,-2.45) {(a)};
	\end{scope}
	\begin{scope}[xshift=5.6cm]
		\pic {celldisk};
		\node at (0,0) {$\varepsilon$};
		\draw[symmetry, -{Stealth[length=2.6mm]}]
			(0:0.6) arc[start angle=0, end angle=180, radius=0.6];
		\node[above] at (90:0.6) {$\pi$};
		\node at (0,-2.45) {(b)};
	\end{scope}
\end{tikzpicture}
\caption{The two symmetries of the cell decomposition of the
$2$-disk. (a)~The reflection, which fixes $a$ and $b$ and
interchanges the two edges: $x\mapsto -y$, $y\mapsto -x$ and
$\varepsilon\mapsto-\varepsilon$ (Proposition \ref{prop:reflection}).
(b)~The rotation by $\pi$, which interchanges the vertices and the
edges: $a\leftrightarrow b$, $x\leftrightarrow y$ and
$\varepsilon\mapsto e^{\ad_{-x}}(\varepsilon)$
(Proposition \ref{prop:rotation}).}
\label{fig:symmetries}
\end{figure}
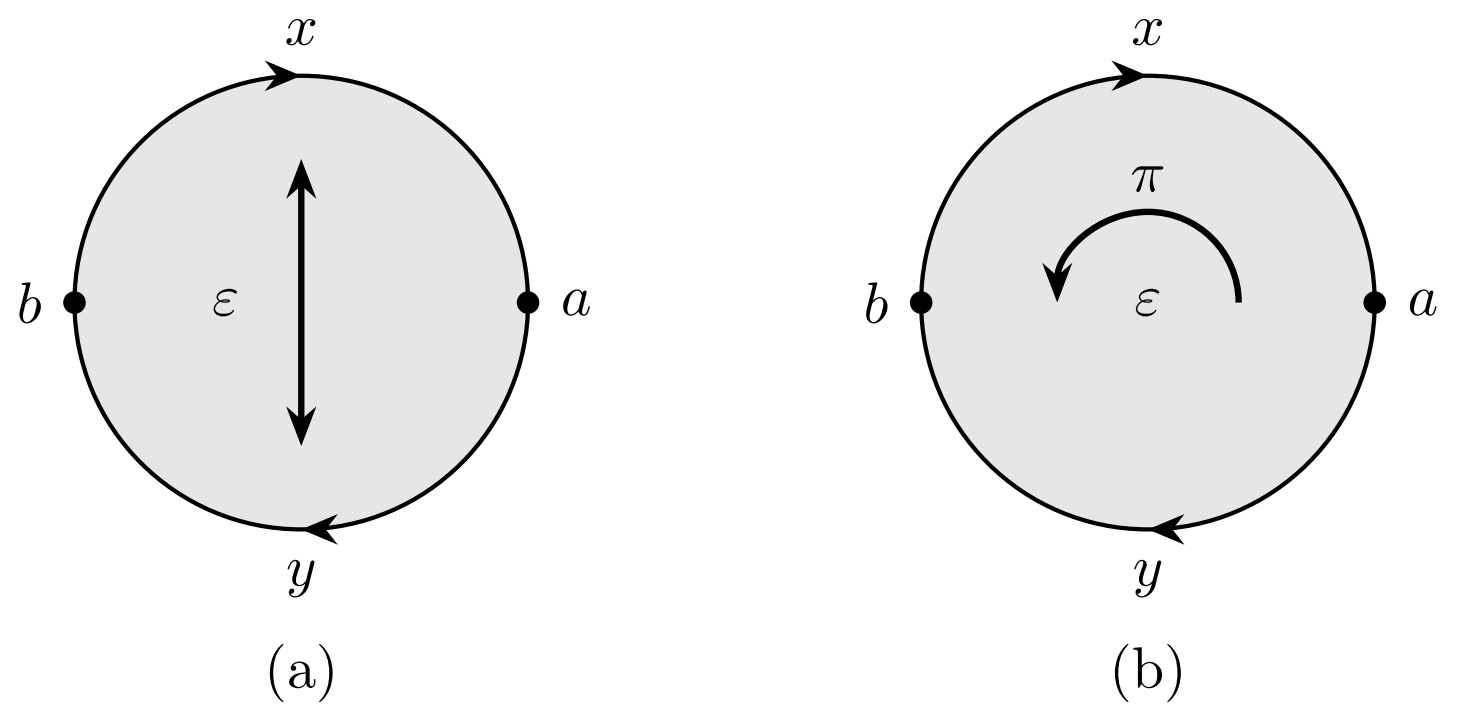

In this section we consider two symmetries of the cell decomposition
of the $2$-disk of Section \ref{sec: Model of the 2-disk}, depicted in
Figure \ref{fig:symmetries}: the reflection that fixes the two
vertices and interchanges the two edges, and the rotation by $\pi$,
which interchanges the vertices and also the edges. Each of them
corresponds to a change of generators of the complete dg Lie algebra
of Theorem \ref{disk} in which the differential is given by the same
formulas. For the reflection, this amounts to the inversion formula of
the BCH product, just as the symmetry of the interval amounts to the
vanishing of the odd Bernoulli numbers
(Section \ref{sec: vanishing of odd Bernoulli}). For the rotation, it
amounts to an instance of the conjugation identity of
Lemma \ref{lem:conjugation}(1), and through this identity the rotation
implies the equality $B_4=B_8$, without using the value of any
Bernoulli number.

Throughout this section,
$\bigl(\widehat{\mathbb{L}}(a,b,x,y,\varepsilon),\partial\bigr)$ is
the complete dg Lie algebra of Theorem \ref{disk}, so that $x$ is a
gauge from $b$ to $a$, $y$ is a gauge from $a$ to $b$, and
$\partial_a\varepsilon=x*y$ (Section \ref{sec: Model of the 2-disk}).
We write $g(s)=s/(e^s-1)=\sum_{n\geq0}B_n s^n/n!$, as in
Remark \ref{rem:associativity}.

Both symmetries are obtained as follows. Let
$a',b',x',y',\varepsilon'$ be elements of
$\widehat{\mathbb{L}}(a,b,x,y,\varepsilon)$ of degrees
$-1,-1,0,0,1$ whose linear parts form a basis of the span of
$a,b,x,y,\varepsilon$, and assume that they satisfy the formulas of
Theorem \ref{disk}:
\begin{equation}\label{eq:diskformulas}
\begin{aligned}
&\partial a'=-\tfrac12[a',a'],\qquad \partial b'=-\tfrac12[b',b'],\qquad \partial_{a'}\varepsilon'=x'* y',\\
&\partial x'=\ad_{x'}(b')+g(\ad_{x'})(b'-a'),\\
&\partial y'=\ad_{y'}(a')+g(\ad_{y'})(a'-b').
\end{aligned}
\end{equation}
Then $a',b',x',y',\varepsilon'$ are free generators of
$\widehat{\mathbb{L}}(a,b,x,y,\varepsilon)$, and in these generators
the differential is literally the one of Theorem \ref{disk}. In other
words, the assignment $a\mapsto a'$, $b\mapsto b'$, $x\mapsto x'$,
$y\mapsto y'$, $\varepsilon\mapsto\varepsilon'$ is an automorphism of
the complete dg Lie algebra
$\bigl(\widehat{\mathbb{L}}(a,b,x,y,\varepsilon),\partial\bigr)$.
Indeed, it defines a morphism of complete graded Lie algebras, which
is an isomorphism because its linear part is, and
(\ref{eq:diskformulas}) says precisely that it commutes with
$\partial$ on the generators.

\subsubsection{The reflection and the inversion formula of the BCH product}
\label{sec: reflection}

\begin{proposition}\label{prop:reflection}
The elements
$$a'=a,\qquad b'=b,\qquad x'=-y,\qquad y'=-x,\qquad \varepsilon'=-\varepsilon$$
satisfy (\ref{eq:diskformulas}). Hence the assignment $a\mapsto a$,
$b\mapsto b$, $x\mapsto-y$, $y\mapsto-x$,
$\varepsilon\mapsto-\varepsilon$ is an automorphism of the complete
dg Lie algebra of Theorem \ref{disk}.
\end{proposition}

\begin{proof}
The formulas for $a'$ and $b'$ are clear. By
Proposition \ref{prop:symmetry}, applied to each of the two edges,
$-y$ is a gauge from $b$ to $a$ and $-x$ is a gauge from $a$ to $b$,
which gives the formulas for $x'$ and $y'$. Finally, the BCH product
satisfies $(-v)* (-u)=-(u* v)$, because
$e^{-v}e^{-u}=(e^{u}e^{v})^{-1}$, and therefore
\[
\partial_{a'}\varepsilon'=\partial_a(-\varepsilon)=-(x* y)=(-y)* (-x)=x'* y'. \qedhere
\]
\end{proof}

Since $\partial_{a'}\varepsilon'=-(x*y)$ and $x'*y'=(-y)*(-x)$, the
proof shows that, once the formulas for $a'$, $b'$, $x'$ and $y'$
hold, the formula for $\varepsilon'$ in (\ref{eq:diskformulas}) is
equivalent to the inversion formula
$$(-y)*(-x)=-(x*y)$$
of the BCH product. Thus the reflection of the disk corresponds to
the inversion formula of the BCH product, just as the symmetry of the
interval corresponds to the conditions $B_1=-\frac{1}{2}$ and
$B_{2k+1}=0$, $k\geq 1$ (Proposition \ref{prop:symmetry}).

\subsubsection{The rotation and the equality \texorpdfstring{$B_4=B_8$}{B4=B8}}
\label{sec: rotation}

For the rotation we need the following two well-known facts, see
\cite[Propositions 4.13 and 4.24]{LieModels}. The first one also
appears in the work of Gadish, Griniasty and Lawrence
\cite[Lemma 1]{GGL}, and in the form $e^{u}ve^{-u}=e^{\ad_u}(v)$ it
is often called Hadamard's lemma.

\begin{lemma}
\label{lem:conjugation}
Let $(L,\partial)$ be a complete dg Lie algebra and let $u,v\in L$ be degree $0$ elements. Then:
\begin{enumerate}
\item[(1)] $e^{\ad_u}(v)=u*v*(-u)$.
\item[(2)] If $a,b\in L$ are Maurer--Cartan elements and $u$ is a gauge from $b$ to $a$, that is,
$\partial u=\ad_u(b)+\sum_{n\geq 0}\frac{B_n}{n!}\ad_u^n(b-a)$, then $$\partial_b\circ e^{\ad_{-u}}=e^{\ad_{-u}}\circ\partial_a.$$
\end{enumerate}
\end{lemma}

\begin{proposition}\label{prop:rotation}
The elements
$$a''=b,\qquad b''=a,\qquad x''=y,\qquad y''=x,\qquad \varepsilon''=e^{\ad_{-x}}(\varepsilon)$$
satisfy (\ref{eq:diskformulas}). Hence the assignment $a\mapsto b$,
$b\mapsto a$, $x\mapsto y$, $y\mapsto x$,
$\varepsilon\mapsto e^{\ad_{-x}}(\varepsilon)$ is an automorphism of
the complete dg Lie algebra of Theorem \ref{disk}.
\end{proposition}

\begin{proof}
The linear part of $\varepsilon''$ is $\varepsilon$. The formulas for
$a''$, $b''$, $x''$ and $y''$ are those of Theorem \ref{disk} for
$b$, $a$, $y$ and $x$, respectively. Finally, since $x$ is a gauge
from $b$ to $a$, parts (2) and (1) of Lemma \ref{lem:conjugation},
applied with $u=x$ and $u=-x$, respectively, give
\begin{align*}
\partial_{a''}\varepsilon''&=\partial_b\,e^{\ad_{-x}}(\varepsilon)=e^{\ad_{-x}}(\partial_a\varepsilon)=e^{\ad_{-x}}(x* y)\\
&=(-x)* x* y* x=y* x=x''* y''. \qedhere
\end{align*}
\end{proof}

Geometrically, $\varepsilon''$ is the $2$-cell $\varepsilon$
transported from $a$ to $b$ by the inverse of the gauge $x$, and it is
attached along the loop $y* x$ based at $b$. Note that the reflection
requires the symmetry of the interval on the edges, whereas for the
rotation the formulas on the edges are automatic. By
Lemma \ref{lem:conjugation}(2) with $u=x$, the formula for $\varepsilon''$ in
(\ref{eq:diskformulas}) is equivalent to the identity
\begin{equation}\label{eq:rotation}
e^{\ad_{-x}}(x*y)=y*x,
\end{equation}
which is the instance of the conjugation identity of
Lemma \ref{lem:conjugation}(1) on which the rotation rests.

Just as the symmetries of the interval are related to the vanishing
of the odd Bernoulli numbers and to the two choices of $B_1$, the
rotation of the disk is related, through identity
(\ref{eq:rotation}), to an identity among Bernoulli numbers. Write
(\ref{eq:rotation}) as $x*y=e^{\ad_x}(y*x)$ in
$\widehat{\mathbb{L}}(x,y)$, and compare the parts linear in
$y$. By Section \ref{sec: generalized Bernoulli and 2-disk}, the part
of $x* y$ linear in $y$ is $g(-\ad_x)(y)$, and by
Remark \ref{rem:associativity}, with $x$ in place of $y$ and $-y$ in
place of $u$, the part of $y*x$ linear in $y$ is
$g(\ad_x)(y)$. Hence $g(-\ad_x)(y)=e^{\ad_x}g(\ad_x)(y)$, and since
the elements $\ad_x^n(y)$, $n\geq 0$, are linearly independent,
$g(-s)=e^sg(s)$. Comparing coefficients,
\begin{equation}\label{eq:shift}
(-1)^nB_n=\sum_{k=0}^{n}\binom{n}{k}B_k,\qquad n\geq 0.
\end{equation}
In the notation of Section \ref{sec: Choice in Bernoulli}, identity
(\ref{eq:shift}) reads $B_n^+=\sum_k\binom{n}{k}B_k^-$.

Together with the vanishing of the odd Bernoulli numbers $B_{2k+1}$,
$k\geq 1$ (Proposition \ref{prop:symmetry}), identity
(\ref{eq:shift}) implies Kaneko's recurrence \cite{Kaneko}
\begin{equation}\label{eq:kaneko}
\sum_{i=0}^{k}\binom{k+1}{i}(k+i+1)B_{k+i}=0,\qquad k\geq 1.
\end{equation}
Indeed, (\ref{eq:shift}) says that the linear map
$\mathfrak b\colon\mathbb{Q}[z]\to\mathbb{Q}$ given by
$\mathfrak b(z^n)=B_n$ is invariant under $z\mapsto -1-z$. Hence
$\mathfrak b$ vanishes on the derivative of
$\bigl(z(z+1)\bigr)^{k+1}$, which changes sign under this
substitution, and this is (\ref{eq:kaneko}) plus the term
$(2k+2)B_{2k+1}=0$.

For $k=3$ and $k=4$, using again the vanishing of the odd Bernoulli
numbers, identity (\ref{eq:kaneko}) reads, after dividing by $4$ and
by $5$,
$$5B_4+7B_6=0\qquad\text{and}\qquad B_4+14B_6+9B_8=0.$$
Subtracting twice the first identity from the second one, we get
$9(B_8-B_4)=0$. Thus the rotation of
the disk, through identity (\ref{eq:rotation}), together with the
symmetry of the interval, implies the equality
$$B_4=B_8.$$
Both numbers are equal to $-\frac{1}{30}$, and this is the only
coincidence between two non-zero Bernoulli numbers with different
indices.\footnote{Since $|B_{2n}|=2(2n)!\,\zeta(2n)/(2\pi)^{2n}$, the
sequence $|B_{2n}|$ decreases for $1\leq n\leq 3$ and increases for
$n\geq 3$. The remaining cases are checked directly.}

\begin{remark}\label{rem:novalues}
The argument above does not use the value of any Bernoulli number.
In particular, the values of $B_2$ and $B_6$ are not needed. This
cannot be achieved with Euler's identity (Example \ref{ex:euler}) or,
more generally, with the identities of Theorem \ref{teo: mainintro}.
Indeed, if a family of numbers satisfies these identities and the
vanishing of $B_{2k+1}$, $k\geq 1$, then so does the family obtained
by the rescaling $B_n\mapsto\lambda^nB_n$,
$\mathcal{B}_w\mapsto\lambda^{|w|-1}\mathcal{B}_w$, for any
$\lambda\in\mathbb{Q}$, where $|w|$ is the length of $w$. However,
the equality $B_4=B_8$ is not preserved. For instance, Euler's
identity for $n=4,6,8$ gives $B_8=-30B_4^2$, and to conclude
$B_4=B_8$ one still needs the value $B_4=-\frac{1}{30}$, or
equivalently $B_2^2=\frac{1}{36}$. The identity $g(-s)=e^sg(s)$
provided by the rotation of the disk is not preserved by the
rescaling, and this is why it yields $B_4=B_8$ directly.
\end{remark}


\bigskip

\noindent\sc{Urtzi Buijs \&  José Manuel Moreno Fernández}\\ 
\noindent\sc{Departamento de Álgebra, Geometría y Topología, \\ Universidad de Málaga, 29080 \\ Málaga, Spain}\\
\noindent\tt{ubuijs@uma.es \& josemoreno@uma.es}\\

\end{document}